\documentclass[10pt,a4paper]{article}
\usepackage{setspace}
\usepackage[english]{babel}		
\usepackage[latin1]{inputenc}
\usepackage{amscd}
\usepackage{amsfonts}
\usepackage{amsgen}
\usepackage{amsmath}
\usepackage{amssymb}
\usepackage{amstext}
\usepackage{amsthm}
\usepackage{graphicx}
\usepackage{latexsym}
\usepackage{makeidx}																
\usepackage{epsfig}
\usepackage{wrapfig}
\usepackage[all,knot,arc,import,poly]{xy}
\usepackage[usenames]{color}
\usepackage[mathscr]{eucal}
\usepackage{indentfirst}
\usepackage[T1]{fontenc}
\usepackage{setspace}
\newtheorem{Def}{Definition}[section]
\newtheorem{Teo}[Def]{Theorem}
\newtheorem{Prop}[Def]{Proposition}
\newtheorem{Lema}[Def]{Lemma}

\newtheorem{Cor}[Def]{Corollary}

\newcommand{\MS}{\mathcal{S}}
\newcommand{\MA}{\mathcal{A}}

\newcommand{\menos}{\backslash}

\newcommand{\setad}{\rightarrow}

\newcommand{\sse}{\Leftrightarrow}

\newcommand{\sigx}{\sigma(x)}
\newcommand{\sigy}{\sigma(y)}
\newcommand{\sigxy}{\sigma(xy)}

\newcommand{\rgmin}{(RG)^-}

\newcommand{\rgm}{(RG)^-}

\newcommand{\uni}{\mathcal U}
\newcommand{\centro}{\mathcal Z}

\newcommand{\inv}{^{-1}}

\newcommand{\dem}{\begin{proof}}
\newcommand{\cqd}{\end{proof}}

\newcommand{\che}{\left\{}
\newcommand{\chd}{\right\}}
\newcommand{\pare}{\left(}
\newcommand{\pard}{\right)}
\newcommand{\lan}{\langle}
\newcommand{\ran}{\rangle}

	\newcommand{\properpagestyle}
	
\makeindex

\title{Anticommutativity of Skew-symmetric Elements under Generalized Oriented Involutions}
\author{Tonucci, L. E. \and Petit Lobão, T. C.}

\begin{document}
\maketitle


\renewcommand{\contentsname}{Índice}

\pagenumbering{arabic}
\setcounter{page}{1}
%
%


\begin{abstract}
 Let $R$ be a ring with $char(R)\neq2$ whose unit group are denoted by $\uni(R)$, $G$ a group, and $RG$ its group ring. Let $*$ be an involution in $G$, $\sigma:G\setad\uni(R)$ be a nontrivial group homomorphism, with $ker\ \sigma=N$, satisfying $xx^*\in N$ for all $x\in G$, and define the generalized oriented involution $\sigma*$ in $RG$ by $\pare \sum_{x\in G}\alpha_xx\pard^{\sigma*}=\sum_{x\in G}\sigx\alpha_xx^*$. An element $\alpha\in RG$ is called skew-symmetric if $\alpha^{\sigma *}=-\alpha$, and the set of all skew-symmetric elements are denoted by $\rgm$. In this paper, we will classify the group rings $RG$ such that $\rgm$ is anticommutative, generalizing, and obtaining as consequence, the main result of \cite{GP13a}.
\vspace{.1cm}\\
{\it Keywords:} Group rings, rings with involution, skew-symmetric elements, generalized oriented involution, polynomial identity.\\
\end{abstract}

\section{Introduction}

Let $RG$ be a group ring of a group $G$ over a commutative ring $R$ with unity. Given $*$, an involution in $G$, we can naturally induce an involution in $RG$, defined by the linear extension of $*$. 

In a ring $R$ with involution $*$, we call skew-symmetric elements those $r\in R$ such that $r^*=-r$ and collect them in the set $\rgm$. In the same way, we collect the symmetric elements, $r\in R$ such that $r^ *=r$, in the set $(RG)^ +$. Many papers classify the group rings in which these sets, defined with the linear extention of a group involution, satisfy a polinomial identity, see \cite{BJPM09,GP13a,JM06,LSS09}, and when a polinomial identity satisfied in these sets could be lifted to the entire group ring, see \cite{GPS09}.

Using an homomorphism $\sigma:G\setad\che\uni(R)\chd$, we can define $\sigma*:RG\setad RG$ mapping $\sum_{x\in G}\alpha_xx\mapsto \sum_{x\in G}\alpha_x\sigx x^*$; easily we can check this map is an involution if, and only if, $x^*x\inv\in ker\sigma$. When $\sigma$ is nontrivial, $\sigma*$ is called a generalized oriented involution.

In the particular case $\sigma(G)=\che\pm1\chd$, $\sigma*$ is simply called an oriented involution and the papers \cite{BP06,BJPM09,GP13b,GP14,CP12}, searching for a generalization, study the identities stated above under this involution. To achieve these results, the authors used the information of the subgroup $N\leq G$, that was given by the similar results under the linearly extended involution. 

Similarly to the fact that results on linear extention give information about $N\leq G$, the special case $\sigma(G)=\che\pm1\chd$ does the same regarding the subgroup $\sigma^{-1}(\che\pm1\chd)$, that is $C$. For instance, the description of $C$ given in \cite{GP14} was relevant to the authors in order to obtain the main theorem in \cite{PT15}, and, although such a description in \cite{CP12} was not actually used in \cite{V13}, the techniques developed in the first article were often used in the second one.


We say that a group $G$ has an unique nontrivial commutator $s$, if $G'=\che 1,s\chd$. This class of group is quite frequent in the study of commutativity and anticommutativity of symmetric and skew-symmetric elements; some examples can be found in \cite{BP06,BJPM09,GP13a,GP14,JM06,PT15}, where the involution, in this case, is given by $x^*=\che x,sx\chd,\ \forall x\in G$. A group $G$ is said to have limited commutativity, LC-group for short, if, given $x,y\in G$ such that $(x,y)=1$, then $x,y$ or $xy\in\centro(G)$, where $\centro(G)$ is the center of $G$. A special involution is naturally endowed in a LC-group $G$ that has an unique nontrivial commutator $s$, namely, $*:G\setad G$ mapping
$$x^*=\che\begin{array}{rl}
x,&\text{if}\ x\in\centro(G)\\
sx,&\text{contrary case;}
\end{array}\right.
$$
obviously this is a particular case of the previous one. We can verify this map defines an involution only in this group class and, in this case, we say that $G$ is a SLC-group with canonical involution $*$. A better description of SLC-groups is given by [JM06, Theorem 2.4]. 

The groups with involution mentioned above shall play an important role in the our main theorem. The goal of this result is to classify the group rings $RG$ in which the set of skew-symmetric elements, over a generalized oriented involution, is anticommutative. In order to do that, we deal with a similar argument used in \cite{PT15}, however we simply use the $N$ desription, instead of that of $C$, to obtain the result; this way we also reach the description of $C$ as a consequence.

The kernel of $\sigma$ will be denoted by $N$ and the center of $G$, by $\centro(G)$. The symmetric elements under $*$ in $G$ are collected in $G_*$; $(x,y)=x\inv y\inv xy$ is the multiplicative commutator; and $x^y=y\inv xy$ is the conjugation of $x$ by $y$. We will denote $N_*=G_*\cap N$ and $A'=\lan (x,y):x,y\in A\ran$. Throughout this work, we will assume that $char(R)\neq2$ and will use this fact without further mention.

\section{Skew-symmetric Elements Anticommute}

In order to prove that $\rgm$ is anticommutative, it is enough to do that to a set of generators of $\rgm$.


For $\alpha\in RG$, write
$$\alpha=\sum_{x\in N_*}\alpha_xx+\sum_{x\in N\menos G_*}\alpha_xx+\sum_{x\in G_*\menos N}\alpha_xx+\sum_{x\in G\menos(G_*\cup N)}\alpha_xx,$$
then
$$\alpha^{\sigma*}=\sum_{x\in N_*}\alpha_x x+\sum_{x\in N\menos G_*}\alpha_x x^*+\sum_{x\in G_*\menos N}\sigx\alpha_x x+\sum_{x\in G\menos(G_*\cup N)}\sigx\alpha_x x^*,$$
so
$$\alpha^{\sigma*}=-\alpha\ \text{if and only if}\che
\begin{array}{ll}
\alpha_x=-\alpha_{x^*}&\text{for}\ x\in N_*\\
\alpha_x=-\alpha_{x^*}&\text{for}\ x\in N\menos G_*\\
\alpha_x=-\sigx\alpha_{x}&\text{for}\ x\in G_*\menos N\\
\alpha_x=-\sigx\alpha_{x^*}&\text{for}\ x\in G\menos(G_*\cup N).\\
\end{array}\right.$$
Thus $(RG)^-$ is spanned over $R$ by elements in the sets,
$$\MA_1=\che\alpha x:x\in G_*\ \text{e}\ \alpha(1+\sigx)=0\chd$$
$$\MA_2=\che x - \sigx x^*:x\in G\menos G_*\chd.$$

Since $N$ is invariant under $*$, we can easily verify that $\sigma*|_{RN}=*|_{RN}$, therefore, the oriented involution acts in $RN$ as an ordinary involution, so we can deduce some information about $N$ using the following result. 

\begin{Teo}[Theorem 2.2, \cite{GP13a}] Let $G$ be a group with involution $*$ and let $R$ be a ring of $char(R)\neq 2$. Then the set $(RG)^-$ is anticommutative if, and only if, either:
\begin{itemize}
 \item [(1)] $G$ is abelian and $*=Id$.
 \item [(2)] $char(R)=4$, $G$ is abelian, and exists $s\in G$ with $s^2=1$ and $x^*\in\che x,sx\chd,\ \forall x\in G$.
 \item [(3)] $char(R)=4$, $G$ is a nonabelian group with a unique nonidentity commutator $s$, and $x^*\in\che x,sx\chd,\ \forall x\in G$.
\end{itemize}
\label{teo.gp13a}
\end{Teo}

\begin{Lema}Suppose that $\rgm$ is anticommutative. If $x\in G_*$ then $\sigx\neq -1$.\label{aa.simet.menos.1}\end{Lema}

\dem If $x\in G_*$ and $\sigx=-1$, then, for all $\alpha\in R$, the equation $\alpha(1+\sigx)=0$ holds, so $\alpha x\in\rgmin$. Thus, taking $\alpha=1$,  and $x\in\rgm$ we obtain that $\alpha x= x$ anticommutes with itself, so $x^2+x^2=0$, which implies $2x^2=0$; a contradiction, for $char(R)\neq2$.\cqd

The proofs of Lemmas \ref{aa.char4} - \ref{aa.2.simet} are quite similar to some Lemmas in \cite{PT15}, so we will just explain the details we can not point to this paper.

\begin{Lema}If $\rgm$ is anticommutative, then $char(R)=4$ or $8$. Furthermore, if $*\neq Id$, then $char(R)=4$, $xx^*=x^*x$ and $x^2\in G_*$, for all $x\notin G_*$.
\label{aa.char4}
\end{Lema}

\dem The proof is analogous to [PT15, Lemmas 3.2 and 3.3].\cqd

   

\begin{Cor}If $\rgm$ is anticommuative, $*\neq Id$, and $x\in G_*$, then there exists a nonzero $\alpha\in R$ such that $\alpha x\in\rgm$.\label{aa.existencia.alpha}\end{Cor}
\dem If $\sigx=1$, then we can easily verify $\alpha=2$ satisfies such a claim, for $char(R)=4$; if $\sigx\neq1$, then $\alpha=1-\sigx\neq0$ satisfies $\alpha(1+\sigx)=(1-\sigx)(1+\sigx)=1-\sigx^2$, and, for $x\in G_*$ and $x^*x\inv\in N$, we obtain $\sigx^2=1$, so $\alpha(1+\sigx)=0$.
\cqd


Given $x,y\notin G_*$, as $\rgm$ is anticommutative, the following equation will play the same role as [PT15, equation (1)];
\begin{equation}
\begin{array}{rcl}
0&=&(x-\sigma(x)x^*)(y-\sigma(y)y^*)+(y-\sigma(y)y^*)(x-\sigma(x)x^*)\\
&=&xy+yx+\sigma(xy)x^*y^*+\sigma(xy)y^*x^*-\sigma(y)xy^*-\sigma(y)y^*x-\sigma(x)yx^*-\sigma(x)x^*y=0.
\end{array}
\label{aa.eq.1}
\end{equation}

\begin{Lema} Suppose that $\MS$ is anticommutative. Given $x,y\in G$, then, $xy=yx$ if, and only if, $x^*y=yx^*$. Moreover, if $x,y\notin G_*$, then $xy=yx=x^*y^*=y^*x^*$, $xy^*=y^*x=x^*y=yx^*$ and $2(1+\sigxy)=2(\sigx+\sigy)=0$.\label{aa.comut}\end{Lema}

\dem The proof is analogous to [PT15, Lemma 3.4].\cqd

\begin{Lema}Suppose that $\rgm$ is anticommutative. If $x \notin G_*$ and $y\notin N\cup G_*$, then, $x^y\in\che x^*,x\chd$.
\label{aa.base.1}\end{Lema}

\dem We can prove this lemma analogously to [PT15, Lemma 3.5] replacing [PT15, equation (1)] with equation (\ref{aa.eq.1}), except in case $xy=x^*y^*$ and $\sigxy=-1$.

Suppose that $x^y\notin \che x,x^*\chd$, $xy=x^*y^*$ and $\sigxy=-1$. Notice that $(xy)^*=y^*x^*=yx$, then $xy+yx=xy-\sigxy (xy)^*\in\rgm$, and
$$\begin{array}{rcl}
   0&=&(xy+yx)(x\inv-\sigma(x\inv) x^{-*})+(x\inv-\sigma(x\inv) x^{-*})(xy+yx)\\
   &=& 2y+xyx\inv+x\inv yx-\sigma(x\inv)(xyx^{-*}+yxx^{-*}+x^{-*}xy+x^{-*}yx),
\end{array}$$
which implies $y\in\che xyx\inv,xyx^{-*},yxx^{-*},x\inv yx, x^{-*}xy,x^{-*}yx\chd$. Since $x\notin G_*$ and $x^y\notin \che x,x^*\chd$, it follows that $yx=x^*y$. As $yx=y^*x^*$, if $yx=x^*y$, we have $y^*x^*=x^*y$, and, applying involution, $xy=y^*x$; which leads to a contradiction by the proof of [PT15, Lemma 3.5].
\cqd

\begin{Lema}Suppose that $\rgm$ is anticommutative. If $x,y\notin G_*$, then it holds:
\begin{itemize}
 \item [(i)] $x^y\in\che x,x^*\chd$. Particularly, $xy=x^*y^*$.
 \item [(ii)] $xy\in G_*\Leftrightarrow xy=yx$.
 \item [(iii)] If $xy\neq yx$, then $1+\sigma(xy)=\sigma(x)+\sigma(y)$. If $xy=yx$, then $2(1+\sigma(xy))=0=2(\sigx+\sigy)$.
\end{itemize}
\label{aa.2.n.simet}
\end{Lema}

\dem Suppose $(x,y)\neq1$. 

If $x,y\in N$, then, $N$ satisfies (C) of Theorem \ref{teo.gp13a}, thus it has a unique nontrivial commutator $s$ and $x^*=sx$, so $x^y=x(x,y)=xs=x^*$. Furthermore, $xy=xssy=x^*y^*$.

If $x,y\notin N$, then, applying Lemma \ref{aa.base.1} to $x$ and $y$, we get $xy=yx^*$; moreover, applying the same lemma to $y$ e $x^*$, we get $yx^*=x^*y^*$, in other words, $xy=x^*y^*$.

If $x\in N$ and $y\notin N$, applying Lemma \ref{aa.base.1} to $x$ and $y$, as well as to $x^*$ and $y^*$, we have $x^y=x^*$ and $x^*y^*=y^*x$. As $x\inv \in N\menos G_*$ and $(x\inv,y)\neq1$, then, $yx\inv =(x\inv)^*y\neq (x\inv)^*y^*=(yx\inv)^*$, in other words, $yx\inv\notin G_*$; since $yx\inv\notin N$, applying the above case to $y$ and $yx\inv$, it follows that $y^{x\inv}=y^{yx\inv}=y^*$, so $xy=y^*x=x^*y^*$. 

If $x\notin N$ and $y\in N$, by the previous case, $yx=xy^*=y^*x^*$, and, applying involution, we find $x^*y^*=yx^*=xy$, so, $x^y=x^*$.

Finnaly, if $xy=yx$, then Lemma \ref{aa.comut}, guarantees (i).

The item (ii) follows as [PT15, Lemma 3.7 (ii)].

If $(x,y)\neq 1$, then applying item (i) to $x$ and $y$, as well as to $y$ and $x$, we have $xy=yx^*=y^*x=x^*y^*$; thus, item (iii) follows as [PT15, Lemma 3.6 (i) and (ii)], using equation (\ref{aa.eq.1}) instead of [PT15, equation (1)].
\cqd

\begin{Lema}Suppose that $\rgm$ is anticommutative. Then, for all $y\in G_*$, $x\notin G_*$ and $\alpha\in R$ such that $\alpha y\in \rgmin$, it follows that:
\begin{itemize}
 \item [(i)] $x^y\in\che x,x^*\chd$.
 \item [(ii)] $xy\in G_*\sse xy\neq yx$.
 \item [(iii)] If $xy\neq yx$, then $\alpha\sigma(x)=\alpha$. If $xy=yx$, then $\alpha=-\alpha$.
 \end{itemize}
\label{aa.um.em.cada}
\end{Lema}

\dem The proof is analogous to [PT15, Lemma 3.9].
\cqd

\begin{Lema} Suppose that $\rgm$ is anticommutative. 
 If $(G_*)'=\che1\chd$, then $G$ is abelian or a SLC-group with canonical involution $*$.
\label{aa.simet.centro.slc}
 \end{Lema}

\dem Let $x\in G_*$ and $y\notin G_*$. If $xy\neq yx$, then, by item (ii) of Lemma \ref{aa.um.em.cada}, $xy\in G_*$; so, by hyphotesis, $(x,xy)=1$, thus $x(xy)=(xy)x$, which implies $xy=yx$, a contradiction; so $xy=yx$ for all $x\in G_*$ and $y\notin G_*$. As, according to the hyphotesis, $(G_*)'=\che1\chd$, then $G_*\subset\centro(G)$.

Let $(RG)_*$ be the spanning over $R$ by the set $G_*\cup\che x+x^*:x\notin G_*\chd$. We will prove that $(RG)_*$ is commutative. Since $G_*\subset\centro(G)$, it is enough to prove that $(x+x^*)$ commutes with $(y+y^*)$, $\forall x,y\notin G_*$. 

If $xy=yx$, we can easily prove that $(x+x^*)(y+y^*)=(y+y^*)(x+x^*)$. Suppose that $xy\neq yx$, so,
$$\begin{array}{rcll}
   (x+x^*)(y+y^*)&=&xy+xy^*+x^*y+x^*y^*\\
   &=&yx^*+y^*x^*+yx+y^*x&\text{(by item (i) of Lemma \ref{aa.um.em.cada})}\\
   &=&(y+y^*)(x+x^*).
  \end{array}$$
Thus $(RG)_*$ is commutative.

Let $*:RG\setad RG$ be the linear extension of $*:G\setad G$ to $RG$, and notice that $(RG)_*$ is the set of symmetric elements under $*$, so applying [JM06, Theorem 2.4], we conclude that $G$ is abelian or a SLC-group with canonical involution $*$.
\cqd

\begin{Lema}Suppose that $\rgm$ is anticommutative. Then, for all $x,y\in G_*$ and $\alpha,\beta\in R$ such that $\alpha x,\beta y\in\rgm$, it holds that:
\begin{itemize}
 \item [(i)] $xy=yx\sse xy\in G_*$.
 \item [(ii)] If $xy\neq yx$, it follows that $\alpha\beta=0$; if $xy=yx$, then $2\alpha\beta=0$.
 \end{itemize}
 \label{aa.2.simet}
\end{Lema}

\dem This proof is analogous to [PT15, Lemma 3.10].
\cqd

\begin{Prop} If $\rgm$ is anticommutative, then the following properties holds:
\begin{itemize}
 \item [(i)] If $*\neq Id$, then exists $s\in\centro(G)\cap G_*$ such that $s^2=1$ and $x^*\in\che x,sx=xs\chd\ \forall x\in G$.
 \item [(ii)] If $G$ is nonabelian, then $G$ has a unique nontrivial commutator $s\in\centro(G)\cap G_*$ (so it has order $2$) and $x^*\in\che x,sx=xs\chd\ \forall x\in G$.
 \end{itemize}
 \label{aa.ucnt}
 \end{Prop}

\dem (i) Suppose that $*\neq Id$, thus $G\menos G_*\neq\emptyset$. Let $x\notin G_*$ and $s=x\inv x^*$. Naturaly $x^*=xs$ and since $(x,x^*)=1$, then $x^*=sx$; furthermore, by Lemma \ref{aa.char4}, $s^2=(x\inv x^*)^2=x^{-2}(x^*)^2=x^{-2}x^2=1$ and $s=xx^{-*}$. Finnaly, if $y\notin G_*$, then, by Lemma \ref{aa.2.n.simet} (i), $xy=x^*y^*$, in other words, $s=xx^{-*}=y\inv y^*$; thus $s$ does not depend on $x$, so $y^*\in\che y,sy=ys\chd\ \forall y\in G$.

Suppose that $(s,x)\neq1$, for some $x\in G$, then $sx\inv sx=(s,x)=s$; thus $s^x=1$, a contradiction; so, $s\in\centro(G)$. Furthermore, if $s\notin G_*$, then $s^*=s^2=1$, also a contradiction, since $1^*=1\neq s$.

(ii) If $G$ is nonabelian, then the identity does not define an involution on $G$, so $Id\neq *$, thus, by (i), if $x\notin G_*$, then $x^*=xs=sx$.

Notice that, if $x\notin G_*$ and $y\in G$ satisfy $(x,y)\neq1$, then, by Lemma \ref{aa.2.n.simet} (i) or Lemma \ref{aa.um.em.cada} (i), we obtain that $x^y=x^*$, so, $(x,y)=x\inv x^y=x\inv x^*=x\inv xs=s$.

If $x,y\in G_*$ and $(x,y)\neq1$, by Lemma \ref{aa.2.simet} (i), $xy\notin G_*$; hence, by item (i), $yx=y^*x^*=(xy)^*=sxy$, in other words $(x,y)=s$.

\cqd

\begin{Teo}Let $R$ be a commutative ring with $char(R)\neq2$, let $G$ be group with involution $*$ and let $\sigma:G\setad\uni(R)$ be a nontrivial orientation which is compatible with $*$ in the sense that $\sigma(xx^*)=1$, for all $x\in G$ and $N=ker\ \sigma$. Then, $\rgm$ is anticommutative if, and only if, the following holds:
\begin{itemize}
 \item [(i)]  One of the these three possibilities holds:
 \begin{itemize}
 \item [(a)] $char(R)=4$ or $8$, $G$ is abelian, and $*=Id_G$.
 \item [(b)] $char(R)=4$, $G$ is abelian, $*\neq Id_G$, and exists $s\in G_*$, such that $s^2=1$ and $x^*\in\che x,xs\chd,\ \forall x\in G$.
 \item [(c)] $char(R)=4$, $G$ has an unique nontrivial commutator $s$ (thus $s\in G_*\cap\centro(G)$) and $x^*\in\che x,xs\chd$.
 \end{itemize}
 \item [(ii)] $\forall x,y\notin G_*$, if $xy\neq yx$, then $1+\sigma(xy)=\sigma(x)+\sigma(y)$; if $xy=yx$, then $2(1+\sigma(xy))=0=2(\sigx+\sigy)$.
 \item [(iii)] $\forall x\notin G_*,\ y\in G_*$, and $\alpha\in R$ with $\alpha y\in \rgmin$. If $xy\neq yx$, then $\alpha\sigma(x)=\alpha$; if $xy=yx$, then $\alpha=-\alpha$.
 \item [(iv)] $\forall x,y\in G_*$ and $\alpha,\beta\in R$ with $\alpha x,\beta y\in \rgmin$. If $xy\neq yx$, then $\alpha\beta=0$; if $xy=yx$, then $2\alpha\beta=0$.
\end{itemize}

Furthermore, if $(G_*)'=\che1\chd$ and $G$ is nonabelian, then $G$ is a SLC-group with canonical involution $*$.
\label{aa.teo.1}
\end{Teo}

\dem Suppose that $\rgm$ is anticommutative. By Lemma \ref{aa.char4}, we conclude that $char(R)=4$ or $8$. If $*=Id_G$, then $G$ is abelian and (a) holds. Suppose $*\neq Id$, then, by Lemma \ref{aa.char4}, $char(R)=4$.

Using Proposition \ref{aa.ucnt} we can easily verify that, if $G$ is abelian, item (i) implies (b), and case contrary, (c) follows from item (ii).


To verify items (ii)-(iv), it is enough apply item (iii), (iii) e (ii), of Lemma \ref{aa.2.n.simet}, \ref{aa.um.em.cada} and \ref{aa.2.simet}, respectively.

The converse could be proved in the same way of converse of [PT15, Theorem 3.15].

Finnaly, if $(G_*)'=\che1\chd$ and $G$ is nonabelian, then Lemma \ref{aa.simet.centro.slc} guaranties that $G$ is a SLC-group with canonical involution $*$.
\cqd

\begin{Cor}[Theorem 2.1 {[GP13]}]Let $R$ be a commutative ring with $char(R)\neq2$, let $G$ be group with involution $*$ and let $\sigma:G\setad\che\pm1\chd$ be a nontrivial orientation which is compatible with $*$ in the sense that $\sigma(xx^*)=1,$ for all $x\in G$ and $N=ker\ \sigma$. Then, $\rgm$ is anticommutative if, and only if, $char(R)=4$ and the following holds:
 \begin{itemize}
  \item [(1)] $G$ is abelian, $*_N=Id_N$ and exists $s\in N_*$ such that $x^*=xs,\ \forall x\notin N$.
  \item [(2)] $G$ is a SLC-group with canonial involution $*$ and $x^*=xs,\ \forall x\notin N$.
 \end{itemize}
 \label{aa.cor.gp13}
\end{Cor}

\dem Suppose $\rgm$ is anticommutative. By Lemma \ref{aa.simet.menos.1},  $G_*\menos N=\emptyset$, so $*\neq Id$. Applying Theorem \ref{aa.teo.1}, we obtain that $G$ is abelian or has an unique nontrivial commutator $s$ and $x^*=sx,\ \forall x\notin G_*$. In particular, as $G_*\menos N=\emptyset$, $s\in N$ and $x^*=sx,\ \forall x\notin N$.

By Theorem \ref{teo.gp13a} we have three possibilities to $N$, namely (A), (B) or (C).

Suppose $G$ is abelian, so (C) does not occur. Let us prove that (B) also does not hold, leading to (A), hence item (1) happens.

If (B) holds, then exists $x\in N$ such that $x^*=sx$. Taking $y\notin N$, we obtain that $y^*=sy$, moreover, $xy\notin N$, which implies $xy\notin G_*$, since $G_*\menos N=\emptyset$. As $x,y,xy\notin G_*$, by Lemma \ref{aa.2.n.simet} (ii), $xy\neq yx$, so $G$ is nonabelian, a contradiction. Thus (1) holds.

Suppose that $G$ is nonabelian. By Theorem \ref{aa.teo.1} and the fact that $G_*\subset N$, in order to find (2), it is enough prove that $\che1\chd=(G_*)'=(N_*)'$. Notice that this condition holds if $N$ is abelian, in other words, if (A) or (B) occur. Suppose that (C) holds and let $x,y\in N_*$; if $xy\neq yx$, by item (ii) of Lemma \ref{aa.2.simet}, $xy\notin G_*\subset N$, a contradiction, for $x,y\in N$. So, (2) holds.

Conversely, notice that items (i) and (ii) of Theorem \ref{aa.teo.1} can be easily verified. Notice that, if $x\in G_*$ and $\alpha x\in\rgmin$ with $\alpha\neq0$, then $\alpha (1+\sigx)=0$, in other words $2\alpha=0$; so, if $G$ is abelian or a SLC-group, then $G_*\subset\centro(G)$, thus, (iii) and (iv) of Theorem \ref{aa.teo.1} hold, so $\rgm$ is anticommutative.
\cqd

 \begin{Prop}Let $R$ be a commutative ring with $char(R)\neq2$, let $G$ be group with involution $*$ and let $\sigma:G\setad\uni(R)$ be a nontrivial orientation which is compatible with $*$ in the sense that $\sigma(xx^*)=1$, for all $x\in G$, with $exp(G/N)=2$ and $N=ker\ \sigma$ satisfying (1) of Theorem \ref{teo.gp13a}. If $G_*\subset N$, then $\rgm$ is anticommutative if, and only if, $char(R)=4$ and one of the following holds:
  \begin{itemize}
 \item [(i)] $G$ is abelian and $\sigma=\che\pm1\chd$.
 \item [(ii)] 
 \begin{itemize}
  \item [(a)] $G$ is a SLC-group with canonical involution $*$ and $G/N\simeq C_2\times C_2$.
  \item [(b)] $\forall x,y\notin G_*$, if $xy\neq yx$, then $1+\sigma(xy)=\sigma(x)+\sigma(y)$; if $xy=yx$, then $2(1+\sigma(xy))=0=2(\sigx+\sigy)$.
  \item [(c)] $2\alpha=0,\ \forall \alpha\in R$ with $\alpha x\in \rgmin$.
  \end{itemize}
 \end{itemize}
  \label{aa.prop.caso.a}
\end{Prop}

\dem Suppose that $\rgm$ is anticommutative. By Theorem \ref{aa.teo.1}, we have $char(R)=4$. Furthermore, since $N$ satisfies (1), $N$ is abelian and $N\subset G_*$, so, if $G_*\subset N$, then, $G_*=N$; as $N$ is abelian, $(G_*)'=\che1\chd$, then, by Theorem \ref{aa.teo.1}, we obtain that $G$ is abelian or a SLC-group with canonical involution $*$, thus $N=G_*\subset\centro(G)$.

Suppose that $G$ is abelian and let $x,y\notin G_*=N$. By item (ii) of Lemma \ref{aa.2.n.simet} we obtain that $xy=yx\sse xy\in N\sse xN=y\inv N=yN$, where the last identity follows due to $exp(G/N)=2$. This way, as $G$ is abelian, any element that is not in $N$, is in the class $xN$, so, $N$ has only two cosets in $G$, which implies that $G/N\simeq C_2$, in other words $\sigma=\che\pm1\chd$. Thus (i) holds.

Suppose now that $G$ is nonabelian, that is, $G$ is a SLC-group. If $G$ is a SLC-group, then $N=G_*=\centro(G)$, so, by [JM06, Theorem 2.4], we have that $G/N=G/\centro(G)\simeq C_2\times C_2$, and (a) holds.

To verify (b), is enough to apply item (ii) of Theorem \ref{aa.teo.1}. Notice that, if $\alpha x\in\rgm$, then $x\in G_*=\centro(G)$, so, applying item (iii) of Theorem \ref{aa.teo.1}, we conclude (c).
 
Conversely, suppose that (i) holds. In this case, $\sigma$ is a classic orienatation, then, by Corolary \ref{aa.cor.gp13}, we have that $\rgm$ is anticommutative.

Finnaly, suppose that (ii) holds. We will prove that item (i)-(iv) of Theorem \ref{aa.teo.1} occur. Item (i) naturally is verified by SLC-groups, and item (ii) is exactly item (b). Observe that, if $G_*=\centro(G)$ and $\alpha x\in\rgmin$, we have that $xy=yx$, $\forall y\in G$. Thus, hyphotesis (c) guaranties (iii) and (iv).
\cqd




%
\singlespacing
\renewcommand{\bibname}{Referências}
\addcontentsline{toc}{chapter}{\bibname}       


\end{document}